\documentclass[11pt,twoside, a4paper,reqno]{amsart}
\pdfoutput=1

\usepackage[shorthands=off,english]{babel}
\usepackage[T1]{fontenc}
\usepackage[utf8]{inputenc}
\usepackage{amsmath,amsfonts,amssymb}
\usepackage{graphicx}
\usepackage{stmaryrd}
\usepackage{float}
\usepackage{caption}
\usepackage{tikz-cd}
\usepackage[]{color}
\usepackage{xcolor}
\usepackage{textcomp}
\usepackage{latexsym}
\usepackage{xspace}
\usepackage{mathrsfs}
\usepackage[babel=false]{csquotes}
\graphicspath{ {images/} }
\MakeOuterQuote{"}
\usepackage{fancyvrb}
\usepackage{appendix}
\usepackage{chngcntr}
\usepackage{etoolbox}
\DeclareMathAlphabet{\mathcal}{OMS}{cmsy}{m}{n}
\definecolor{mymauve}{rgb}{0.4,0,1}
\usepackage[hidelinks,hypertexnames=false]{hyperref}
\hypersetup{
  colorlinks   = true, 
  urlcolor     = blue, 
  linkcolor    = blue, 
  citecolor   = mymauve 
}
\usepackage{amsthm}
\usepackage[english]{cleveref}
\usepackage{enumitem}

\usepackage[margin=1in]{geometry}
\theoremstyle{plain}
\newtheorem{theorem}{Theorem}[section]{}{}
\newtheorem*{theorem*}{Theorem}
\newtheorem{maintheorem}{Theorem}

\newtheorem{lemma}[theorem]{Lemma}

\newtheorem{proposition}[theorem]{Proposition}

\theoremstyle{definition}

\theoremstyle{remark}

\newtheoremstyle{TheoremNum}
    {8pt}{8pt}                            
    {\itshape}                            
    {}                                    
    {\bfseries}                           
    {.}                                   
    { 5pt plus 1pt minus 1pt }            
    {\thmname{#1}\thmnote{ \bfseries #3}} 
\theoremstyle{TheoremNum}

\makeatletter
\def\cref@thmoptarg[#1]#2#3#4{%
    \ifhmode\unskip\unskip\par\fi%
    \normalfont%
    \trivlist%
    \let\thmheadnl\relax%
    \let\thm@swap\@gobble%
    \thm@notefont{\fontseries\mddefault\upshape}%
    \thm@headpunct{.}
    \thm@headsep 5\p@ plus\p@ minus\p@\relax%
    \thm@space@setup%
    #2
    \@topsep \thm@preskip               
    \@topsepadd \thm@postskip           
    \def\@tempa{#3}\ifx\@empty\@tempa%
      \def\@tempa{\@oparg{\@begintheorem{#4}{}}[]}%
    \else%
      \refstepcounter[#1]{#3}
      \@namedef{cref@#3@alias}{#1}
      \def\@tempa{\@oparg{\@begintheorem{#4}{\csname the#3\endcsname}}[]}%
    \fi%
    \@tempa}%
\makeatother

\crefname{corollary}{Corollary}{Corollaries}
\crefname{conjecture}{Conjecture}{Conjectures}
\crefname{theorem}{Theorem}{Theorems}
\crefname{maintheorem}{Theorem}{Theorems}
\crefname{proposition}{Proposition}{Propositions}
\crefname{lemma}{Lemma}{Lemmas}
\crefname{example}{Example}{Examples}
\crefname{definition}{Definition}{Definitions}
\crefname{remark}{Remark}{Remarks}
\crefname{figure}{Figure}{Figures}
\crefname{table}{Table}{Tables}
\crefname{section}{Section}{Sections}
\crefname{appendix}{Appendix}{Appendices}

\newcommand{\tq}{\mathrel{{\ensuremath{\: : \: }}}}
\def\wt{\widetilde}
\def\wh{\widehat}

\def\rightAction{\curvearrowleft}

\def\unsplitExtension{\cdot}

\newcommand{\centralizer}{{\mathbf{C}}}
\newcommand{\centerOfGroup}{{\mathbf{Z}}}

\def\GL{\mathrm{GL}}
\def\PSL{\mathrm{PSL}}
\def\SL{\mathrm{SL}}
\def\lieGroup{\mathbb{G}}
\def\U{\mathbf{U}}

\def\Ad{\mathrm{Ad}}
\newcommand{\Spec}{\Lambda}

\def\C{\mathbb{C}}
\def\F{\mathbb{F}}

\def\XOS{X_1^{OS}}
\def\XOSk{X_1^{OS+k}}
\def\inclusionBrown{i}

\def\MM{\mathcal{M}}
\def\MMbar{\overline{\mathcal{M}}}

\def\HH{\mathcal{H}}

\def\rhoG{{\rho_0}}

\def\identityLieGroup{\mathbf{1}}
\def\identityMM{\identityLieGroup}

\def\onlyEdgeInEMinusT{{\tilde{\eta}}}

\title{Group actions on contractible $2$-complexes II}
\author[K.I. Piterman]{Kevin Iv\'an Piterman$^{1,2}$}
\author[I. Sadofschi Costa]{Iv\'an Sadofschi Costa$^1$}

\address{$^1$ Departamento  de Matem\'atica - IMAS \\
 FCEyN, Universidad de Buenos Aires. Buenos Aires, Argentina.}

\address{$^2$ Department of Mathematics and Data Science \\
Vrije Universiteit Brussel, Pleinlaan 2, 1050 Brussel.}

\email{kevin.piterman@vub.be}
\email{isadofschi@dm.uba.ar}



\thanks{Researchers of CONICET. 
Kevin Iv\'an Piterman was partially supported by grants PIP 11220170100357, PICT 2017-2997, and UBACYT 20020160100081BA, and by an Oberwolfach Leibniz Fellowship.
Iv\'an Sadofschi Costa was partially supported by grants PICT-2017-2806, PIP 11220170100357CO and UBACyT 20020160100081BA}

\begin{document}
\begin{abstract}
 In this second part we prove that,
 if $G$ is one of the groups $\PSL_2(q)$ with $q>5$ and $q\equiv 5\pmod {24}$ or $q\equiv 13 \pmod{24}$,
 then the fundamental group of every acyclic $2$-dimensional, fixed point free and finite $G$-complex admits a nontrivial representation in a unitary group $\U(m)$.
 This completes the proof of the following result: every action of a finite group on a finite and contractible $2$-complex has a fixed point.
\end{abstract}

\maketitle

\setcounter{tocdepth}{1}
\tableofcontents

\section{Introduction}

In this second part we prove the following:

\addtocounter{maintheorem}{2}
\begin{maintheorem}\label{thmC}
 Let $G$ be one of the groups $\PSL_2(q)$ with $q>5$ and $q\equiv 5\pmod {24}$ or $q\equiv 13 \pmod{24}$.
 Then the fundamental group of every $2$-dimensional, fixed point free, finite and acyclic $G$-complex admits a nontrivial representation in a unitary group $\U(m)$.
\end{maintheorem}
 This completes the proof of the following result: every action of a finite group $G$ on a finite and contractible $2$-complex $X$ has a fixed point.

The groups $G$ considered in \cite[Theorem B]{part1} share a key property: they admit a nontrivial representation $\rhoG$ which restricts to an irreducible representation on the Borel subgroup.
The moduli $\MM_k$ of representations of $\Gamma_k = \XOSk(G)\unsplitExtension G$ constructed in the proof of \cite[Theorem B]{part1} is built from a representation with this property.
When $q\equiv 1\pmod 4$ no nontrivial representation of $\PSL_2(q)$ restricts to an irreducible representation on the Borel subgroup.
To prove \cref{thmC}, we circumvent this difficulty by instead considering the action of $\wh{G}=\SL_2(q)$ on $\XOSk(G)$.

\bigskip

\textbf{Acknowledgements.}
This work was partially done during a stay of the first author at The Mathematisches Forschungsinstitut Oberwolfach.
He is very grateful to the MFO for their hospitality and support.

\section{More representation theory}

  We denote the set of eigenvalues of a square matrix $M$ by $\Spec(M)$.
 \begin{lemma}\label{LowerBoundForDimensionOfIntersectionOfCentralizersCyclicSubgroups}
    Let $G$ be a finite group, $g_1,g_2\in G$ and $H_i = \langle g_i\rangle$.
    Let $\rho\colon G\to \U(m)$ be a unitary representation and let $n_i=\# \Spec(\rho(g_i))$.
    Then there are matrices $A_1,A_2\in \U(m)$ such that $A_i\rho(g_i)A_i^{-1}$ is diagonal (for $i=1,2$), $A_1^{-1}A_2$ commutes with $\centralizer_{\U(m)}( \rho(G) )$ and
    \[\dim\left( \left(A_1\centralizer_{\U(m)}\left(\rho(H_1)\right)A_{1}^{-1}\right)\cap \left(A_2\centralizer_{\U(m)}\left(\rho(H_2)\right)A_{2}^{-1}\right)\right) \geq \frac{m^2}{n_1n_2}.\]
  
  \begin{proof}
    We can take $T\in \U(n)$ and irreducible representations $\rho_j\colon G\to \U(m_j)$ with $j=1,\ldots,k$ such that
    $T\rho T^{-1} = \rho_1\oplus \cdots \oplus \rho_k$.
    Moreover, we can do this so that whenever $\rho_j$ and $\rho_{j'}$ are isomorphic we have $\rho_j = \rho_{j'}$.
    For each $i=1,2$ we take matrices $D_{i,1},\ldots,D_{i,k}$ with $D_{i,j}\in \U(m_j)$ such that
    \[D_{i,j} \rho_j(g_i) D_{i,j}^{-1}\]
    is diagonal.
    We choose the $D_{i,j}$ so that $\rho_j=\rho_{j'}$ implies $D_{i,j}=D_{i,j'}$.
    Let $D_i = D_{i,1}\oplus\cdots \oplus D_{i,k}$.
    Then, by \cite[Proposition 7.3 and Remark 7.4]{part1}, $D_i$ commutes with $\centralizer_{\U(m)}( T\rho(G)T^{-1})$ and letting $A_i=D_iT$, we have that $A_i\rho(g_i)A_i^{-1}$ is diagonal and $A_1^{-1}A_2$ commutes with $\centralizer_{\U(m)}(\rho(G))$.
    Now for $\lambda_1\in \Spec(\rho(g_1))$ and $\lambda_2\in \Spec(\rho(g_2))$ we define
    \[n(\lambda_1,\lambda_2)=\# \{ j \tq 1\leq j \leq m \text{ and } (A_1\rho(g_1)A_1^{-1})_{j,j} = \lambda_1 \text{ and } (A_2 \rho(g_2)A_2^{-1})_{j,j} = \lambda_2\}.\]
    Note that
    \[\left(A_1\centralizer_{\U(m)}\left(\rho(H_1)\right)A_{1}^{-1}\right)\cap \left(A_2\centralizer_{\U(m)}\left(\rho(H_2)\right)A_{2}^{-1}\right)\]
    has a subgroup isomorphic to
    \[\prod_{\lambda_1\in \Spec(\rho(g_1)),\lambda_2\in \Spec(\rho(g_2))} \U(n(\lambda_1,\lambda_2))\]
    and therefore its dimension is at least $\sum_{\lambda_1\in \Spec(\rho(g_1)), \lambda_2\in\Spec(\rho(g_2))} n(\lambda_1,\lambda_2)^2$.
    The AM-QM inequality gives
    \[\frac{m}{n_1n_2} =\frac{\displaystyle\sum_{\lambda_1\in \Spec(\rho(g_1)),\lambda_2\in\Spec(\rho(g_2))} n(\lambda_1,\lambda_2)}{n_1n_2}
      \leq \sqrt{\frac{\displaystyle\sum_{\lambda_1\in \Spec(\rho(g_1)),\lambda_2\in\Spec(\rho(g_2))} n(\lambda_1,\lambda_2)^2}{n_1n_2}},\]
    and we obtain the desired inequality.
  \end{proof}
 \end{lemma}

\section{The \texorpdfstring{$\wh{G}$}{hat G}-graph \texorpdfstring{$\XOS(G)$}{XOS1(G)}}

Let $G=\PSL_2(q)$ with $q\equiv 5\pmod{24}$ or $q\equiv13\pmod{24}$.
We consider a construction of $\XOS(G)$ as in \cite[Proposition 3.10]{part1}.
Recall that for any $k\geq 0$, we can also consider the $G$-graph $\XOSk(G)$ obtained from $\XOS(G)$ by attaching $k$ free orbits of $1$-cells.

Let $\wh{G} = \SL_2(q)$, so that $\centerOfGroup(\wh{G}) = \{1,-1\}$ and $\wh{G}/\centerOfGroup(\wh{G}) = G$.
We consider the action of $\wh{G}=\SL_2(q)$ on $\XOSk(G)$ defined using the projection $\pi\colon \wh{G}\to G$.
The stabilizer of a vertex (resp. edge) for the action of $\wh{G}$ is a central extension, by $\centerOfGroup(\wh{G})$, of the stabilizer for the action of $G$.
Then the $\wh{G}$-orbits are connected as in \cref{FigureXOSG}.
The group $B$ denotes the Borel subgroup of $\wh{G}$ and $Q_8$ denotes the quaternion group.

\begin{figure}[H]
   \begin{minipage}{0.48\textwidth}
     \centering
      \begin{tikzcd}
      B \arrow[-]{d}[swap]{C_{q-1}} \arrow[-]{r}[]{C_6}& \SL_2(3) \arrow[-]{dl}[]{Q_8} \\
      2D_{q-1} \arrow[-]{r}[swap]{C_4} & 2D_{q+1}
      \end{tikzcd}

      \bigskip
      
      $\wh{G}=\SL_2(q)$, $q\equiv 13\pmod {24}$.
   \end{minipage}
  \begin{minipage}{0.48\textwidth}
    \centering
    \begin{tikzcd}
    B \arrow[-]{d}[swap]{C_{q-1}}\arrow[-]{r}[]{C_4}  & 2D_{q+1} \arrow[-]{d}[]{C_6} \\
    2D_{q-1} \arrow[-]{r}[swap]{Q_8} & \SL_2(3)
    \end{tikzcd}

     \bigskip
    
    $\wh{G}=\SL_2(q)$, $q\equiv 5\pmod {24}$.
   \end{minipage}
   \caption{The $\wh{G}$-graph $\XOS(G)$.}
   \label{FigureXOSG}
\end{figure}
Below we describe the choices needed, in each case, to  apply Brown's result to the action of $\wh{G}$ on $\XOSk(G)$ (see \cite[Theorem 5.1]{part1}).
The stabilizers for the representatives are given in \cref{TableStabilizersVerticesOSGraph,TableStabilizersEdgesOSGraph}.
\begin{itemize}
\item For $\wh{G}=\SL_2(q)$ with $q\equiv 13\pmod {24}$
we take $V=\{v_0,v_1,v_2,v_3\}$, $E=\{\eta_0,\eta_1,\eta_2,\allowbreak\eta_3,\allowbreak\eta'_1,\ldots,\allowbreak\eta'_k\}$, $T=\{\eta_0,\eta_1,\eta_2\}$, with
$v_0\xrightarrow{\eta_0}v_1$, $v_1\xrightarrow{\eta_1}v_2$,  $v_1\xrightarrow{\eta_2}v_3$, $v_3\xrightarrow{\eta_3} g_{\eta_3} v_0$ and $v_0\xrightarrow{\eta'_i} v_0$ for $i=1,\ldots, k$.
\item For $\wh{G}=\SL_2(q)$ with $q\equiv 5\pmod {24}$
we take $V=\{v_0,v_1,v_2,v_3\}$, $E=\{\eta_0,\eta_1,\eta_2,\allowbreak\eta_3,\allowbreak\eta'_1,\ldots,\eta'_k\}$, $T=\{\eta_0,\eta_2,\eta_3\}$, with
$v_0\xrightarrow{\eta_0}v_1$, $v_1\xrightarrow{\eta_2}  v_3$, $v_3\xrightarrow{\eta_3} v_2$, $v_2\xrightarrow{\eta_1}g_{\eta_1}v_0$ and $v_0\xrightarrow{\eta'_i} v_0$ for $i=1,\ldots, k$.
\end{itemize}
\renewcommand{\arraystretch}{1.25}
\begin{table}[h]
 \centering
\begin{tabular}{|c|c|c|c|c|c|}
  \hline
  ${\wh{G}}$ & ${q}$ & $\wh{G}_{v_0}$ & $\wh{G}_{v_1}$ & $\wh{G}_{v_2}$ & $\wh{G}_{v_3}$  \\
  \hline
  $\SL_2(q)$ & $q$ odd & $B=\F_q\rtimes C_{q-1}$ & $2D_{q-1}$ & $2D_{q+1}$ & $\SL_2(3)$   \\
  \hline  
 \end{tabular}
 \caption{Stabilizers of vertices for the $\wh{G}$-graph $\XOSk(G)$.}
\label{TableStabilizersVerticesOSGraph}
\end{table}

\begin{table}[h]
 \centering
\begin{tabular}{|c|c|c|c|c|c|c|}
  \hline
  $\wh{G}$ & ${q}$  & $\wh{G}_{\eta_0}$ & $\wh{G}_{\eta_1}$ & $\wh{G}_{\eta_2}$ & $\wh{G}_{\eta_3}$ &$\wh{G}_{\eta'_i}$ \\
  \hline
  $\SL_2(q)$ & ${q}$ odd & $C_{q-1}$ & $C_4$ & $Q_8$ & $C_6$  & $\centerOfGroup({\wh{G}})$ \\
  \hline  
 \end{tabular}
 \caption{Stabilizers of edges for the $\wh{G}$-graph $\XOSk(G)$.}
\label{TableStabilizersEdgesOSGraph}
\end{table}
\renewcommand{\arraystretch}{1}
In what follows, $\Gamma_k$ denotes the group obtained by applying Brown's theorem to the action of $\wh{G}$ on $\XOSk(G)$ with these choices.

\section{Representations and centralizers}

In this section we fix a suitable irreducible representation $\rhoG\colon \wh{G}\to \lieGroup$. Recall that $(x)$ denotes the conjugacy class of $x$.

\begin{proposition}\label{rep1mod4}
  Let $\wh{G} = \SL_2(q)$ with $q\equiv 1 \pmod{4}$.
  Then there are elements $a,b,c,d \in G$ with orders $|a|=q-1$, $|b|=q+1$, and $|c|=|d|=p$, such that the following hold:
  \begin{enumerate}[label=(\roman*)]
    \item There are exactly $q+4$ conjugacy classes in $G$:
    $(1)$; $(-1)$; $(a^i)$ for $1\leq i \leq (q-3)/2$; $(b^j)$ for $1\leq j \leq (q-1)/2$; $(c)$; $(-c)$; $(d)$; and $(-d)$.
    \item The elements in a Borel subgroup $B$ are the following:
      $1$; $-1$;
      $2q$ elements in each class $(a^l)$, for $l=1,\ldots,\frac{q-3}{2}$;
      and $\frac{q-1}{2}$ elements in each of the classes $(c)$, $(-c)$, $(d)$ and $(-d)$.
    \item The elements in a $C_4$ subgroup are the following: $1$; $-1$; and two elements in $(a^{\frac{q-1}{4}})$.
    \item The elements in a $2D_{q+1}$ subgroup are the following:
      $1$; $-1$;
      $q+1$ elements in $(a^{\frac{q-1}{4}})$;
      and two elements in each class $(b^m)$ for $m=1,\ldots,\frac{q-1}{2}$.
    \item There is an irreducible character $\chi$ given by
    \begin{center}
      \begin{tabular}{c|cccccccc}
              & $1$       & $-1$       & $(a^i)$ & $(b^j)$      & $(c)$              & $(d)$             \\ \hline 
        $\chi$ & $(q-1)/2$ & $-(q-1)/2$ & $0$     & $(-1)^{j+1}$ & $(-1+\sqrt{q})/2$  & $(-1-\sqrt{q})/2$ \\
      \end{tabular}
  \end{center}
  and we have $\chi(-c)=-\chi(c)$, $\chi(-d)=-\chi(d)$.
  \end{enumerate}
  
  \begin{proof}
    The description of the conjugacy classes and the character table can be found in \cite[Theorem 38.1]{DornhoffPartA}.
    Moreover (ii) is proved in \cite[p. 231]{DornhoffPartA}.
    Now (iii) follows from $(a^{\frac{q-1}{4}})$ being the only conjugacy class of elements of order $4$.
    To prove (iv) note that each order $2$ element in $D_{q+1}\leq \PSL_2(q)$ lifts to two order $4$ elements in $2D_{q+1}$
    and that the order $(q+1)/2$ cyclic subgroup of $D_{q+1}$ lifts to an order $q+1$ cyclic subgroup.
  \end{proof}
\end{proposition}

\begin{proposition}\label{descriptionCentralizersOneModFour}
  Let $\wh{G}=\SL_2(q)$ with $q \equiv 5\pmod 8$ and let $\lieGroup=\U\left((q-1)/2\right)$.
  There is an irreducible representation $\rhoG\colon \wh{G}\to \lieGroup$ satisfying the following properties:
  \begin{enumerate}[label=(\roman*)]
    \item The restriction of $\rhoG$ to the Borel subgroup $\wh{G}_{v_0}$ is irreducible.
    \item The centralizer $\centralizer_{\lieGroup}(\rhoG(\wh{G}_{\eta_1}))$ has dimension $\frac{(q-1)^2}{8}$.
    \item The centralizer $\centralizer_{\lieGroup}(\rhoG(\wh{G}_{v_2}))$ has dimension $\frac{q-1}{4}$.
  \end{enumerate}
\begin{proof}
  We take $\rhoG$ realizing the degree $(q-1)/2$ irreducible character $\chi$ in part (v) of \cref{rep1mod4}.
  By \cite[Theorem 7.1]{part1}, we can take $\rhoG$ to be unitary.
  By \cite[Lemma 7.5]{part1} we can prove parts (i) to (iii) by restricting $\chi$ and computing the norm using parts (ii) to (iv) of \cref{rep1mod4}.
\end{proof}
\end{proposition}

\section{The proof of \texorpdfstring{\cref{thmC}}{Theorem C}}
For each of the groups $G$ in \cref{thmC},
we consider a closed edge path $\xi$ in $\XOS(G)$ such that attaching a free $G$-orbit of $2$-cells along this path gives an acyclic $2$-complex.
We define $x_0=\inclusionBrown(\xi)$, where $\inclusionBrown\colon \pi_1(\XOS(G),v_0)\to \Gamma_0$ is the inclusion given by Brown's theorem.
We set $x_i=x_{\eta'_i}$ for $i=1,\ldots,k$.
Let $\onlyEdgeInEMinusT$ be the unique edge of $\XOS(G)$ which lies in $E-T$.
We define $y_0=x_{\onlyEdgeInEMinusT}$ and $y_i=x_{\eta'_i}$ for $i=1,\ldots,k$.\medskip

Let $\MM_k$ be the moduli of representations of $\Gamma_k$ obtained from the representation $\rhoG\colon \wh{G}\to \U(m)$ of \cref{descriptionCentralizersOneModFour} using \cite[Theorem 5.2]{part1}.
Let $\MMbar_k$ be the corresponding quotient obtained using \cite[Theorem 5.3]{part1}.
Note that the equalities $\MM_k = \MM_0 \times \lieGroup^k$ and $\MMbar_k=\MMbar_0\times \lieGroup^k$ still hold, because $\rhoG(-1)\in\centerOfGroup(\lieGroup)$.

In what follows we consider the induced maps $X_i(\tau)= \rho_\tau(x_i)$, $Y_i(\tau)=\rho_\tau(y_i)$.

\begin{proof}[Proof of \cref{thmC}]
By \cite[Corollary 5.4]{part1}, $\MM_k$ and $\MMbar_k$ are connected and orientable.
Moreover $\dim \MMbar_k =\dim \lieGroup^{k+1}$ (the same proof of \cite[Proposition 8.3]{part1} works in this case).
By \cref{IdentityRegularPoint4k+1} below, $\identityMM$ is a regular point of $X_0$.
By \cref{degreeY05mod24,degreeY013mod24} below, $\overline{Y}_0\colon \MMbar_0\to \lieGroup$ has degree $0$.
The rest of the proof is identical to that of \cite[Theorem B]{part1}.
\end{proof}

\section{The differential of \texorpdfstring{$X_0$}{X0} at \texorpdfstring{$\identityMM$}{1}}

The following lemma extends \cite[Lemma 9.1]{part1} for the action of $\wh{G}$ on $\XOS(G)$.
We denote $z=-1\in\centerOfGroup(\wh{G})$ to avoid confusion with $-1\in \C\subseteq \C[\wh{G}]$.

\begin{lemma}
\label{complementLinearCombinationFixedPoints}
Let $G$ be one of the groups in \cref{thmC}.
Let $E$ be a set of representatives of the orbits of edges in $\XOS(G)$.
Let $X$ be an acyclic $2$-complex obtained from $\XOS(G)$ by attaching a free orbit of $2$-cells along the $G$-orbit of a closed edge path
$\xi = (a_1e_1^{\varepsilon_1}, \ldots, a_n e_n^{\varepsilon_n})$ with $e_i\in E$, $a_i\in \wh{G}$ and $\varepsilon_i\in \{-1,1\}$.
Then it is possible to choose, for each $e\in E$, an element $x_e\in \C[\wh{G}]$ and an element $\delta\in \C[\wh{G}]$ so that 
\[1 = (1-z)\delta + \sum_{i=1}^n \varepsilon_i a_i N(\wh{G}_{e_i}) x_{e_i}.\]
Therefore, for any complex representation $V$ of $\wh{G}$ we have
$V = (1-z) V + \sum_{e\in E} s_e V^{\wh{G}_e}$, 
where $s_e = \sum_{i\in I_e} \varepsilon_i a_i$ and $I_e=\{i\tq e_i = e\}$.

\begin{proof}
Consider the ring homomorphism $\pi\colon \C[\wh{G}]\to \C[G]$.
By \cite[Lemma 9.1]{part1}, there are elements $\wt{x}_e\in \C[\wh{G}]$ such that
$1 = \sum_{i=1}^n \varepsilon_i \pi(a_i) N(G_{e_i}) \pi(\wt{x}_{e_i})$.
Let $x_e=\frac{1}{2}\wt{x}_e$.
Note that $\pi(N(\wh{G}_e))=2\cdot N(G_e)$ and then $\pi( \sum_{i=1}^n \varepsilon_i a_i N(\wh{G}_{e_i})x_{e_i} ) = 1$.
Therefore, since the kernel of $\pi$ is the ideal generated by $1-z$, there is an element $\delta\in \C[\wh{G}]$ such that $ 1 = (1-z)\delta + \sum_{i=1}^n \varepsilon_i a_i N(\wh{G}_{e_i}) x_{e_i}$.
\end{proof}
\end{lemma}

\begin{proposition}\label{zActsTrivially}
 The representation $\rhoG$ satisfies $\Ad \circ \rhoG(z)=1$, where $\Ad\colon \lieGroup \to \GL(T_\identityLieGroup \lieGroup )$ is the adjoint representation.
 \begin{proof}
  This is immediate, for $\Ad(g)$ is the differential of the map $x\mapsto gxg^{-1}$ and $\rhoG(z)=-1$ is central.
 \end{proof}
\end{proposition}

\begin{lemma}\label{IdentityRegularPoint4k+1}
For each of the groups in \cref{thmC}, $\identityMM$ is a regular point of $X_0\colon \MM_0\to \lieGroup$.
\begin{proof}
 Consider the adjoint representation
 $\Ad\circ \rhoG \colon \wh{G}\to \GL( T_\identityLieGroup \lieGroup )$
 which is given by
 $g\cdot v = d_{\rhoG(g)^{-1}} L_{\rhoG(g)}\circ d_\identityLieGroup R_{\rhoG(g)^{-1}} ( v )$.
 By \cite[Proposition 4.4]{part1}, we have $T_\identityLieGroup \centralizer_\lieGroup( \rhoG( H ) ) = (T_\identityLieGroup \lieGroup)^H$.
 By \cref{zActsTrivially} we have $(1-z)\cdot T_\identityLieGroup \lieGroup = 0$ and then, from \cref{complementLinearCombinationFixedPoints}, we conclude that
 $T_\identityLieGroup \lieGroup= \sum_{e\in E} s_e \cdot T_\identityLieGroup \centralizer_{\lieGroup}(\rhoG(\wh{G}_e))$.
 Now the result follows from \cite[Theorem 5.7]{part1}.
\end{proof}
\end{lemma}

\section{The degree of \texorpdfstring{$\overline{Y}_0$}{bar Y0}}
We now prove the degree of $\overline{Y}_0$ is $0$ for each of the groups in \cref{thmC}.
\cref{TableY0}  gives the value of $Y_0$ in the different cases that we consider.
We modify the argument in \cite[Section 10]{part1}.
For $q\equiv 13\pmod {24}$ the argument is easier, but the case of $q\equiv 5\pmod {24}$ requires a more careful analysis.

\renewcommand{\arraystretch}{1.25}
\begin{table}[H]
 \centering
  \begin{tabular}{|c|c|}
  \hline
  ${q}$ & $Y_0(\tau)$ \\
  \hline
  $q\equiv 13 \pmod {24}$  & $\tau_{\eta_0}^{-1}\tau_{\eta_{2}}^{-1}\tau_{\eta_3}^{-1}\rhoG(g_{\eta_3})$\\
  \hline
  $q\equiv 5 \pmod{24}$ & $\tau_{\eta_0}^{-1}\tau_{\eta_{2}}^{-1}\tau_{\eta_{3}}^{-1}\tau_{\eta_1}^{-1}\rhoG(g_{\eta_1})$ \\
  \hline
  \end{tabular}
  \caption{The map $Y_0$, for each of the cases that we consider.}
  \label{TableY0}
 \end{table}
 \renewcommand{\arraystretch}{1}

 \begin{proposition}\label{degreeY013mod24}
  If $q\equiv 13\pmod {24}$ then the map $\overline{Y}_0\colon \MMbar_0\to \lieGroup$ has degree $0$.
  \begin{proof}
   Consider the manifold $M=\centralizer_{\lieGroup}(\rhoG(\wh{G}_{\eta_0})) \times \centralizer_{\lieGroup}(\rhoG(\wh{G}_{\eta_2})) \times \centralizer_{\lieGroup}(\rhoG(\wh{G}_{\eta_3}))$, the group $H=\centralizer_{\lieGroup}(\rhoG(\wh{G}_{v_1}))\times \centralizer_{\lieGroup}(\rhoG(\wh{G}_{v_3}))$ and the free right action $M\rightAction H$ given by
   \[(\tau_{\eta_0},\tau_{\eta_{2}},\tau_{\eta_3})\cdot (\alpha_{v_1},\alpha_{v_3}) = (\alpha_{v_1}^{-1}\tau_{\eta_0}, \alpha_{v_3}^{-1}\tau_{\eta_{2}}\alpha_{v_1},\tau_{\eta_3}\alpha_{v_3}).\]
   Then by \cref{descriptionCentralizersOneModFour}
   \begin{align*}
    \dim M/H &=\dim M -\dim H\\
    &= \dim \MM_0 -\dim \HH +\dim \centralizer_\lieGroup(\rhoG(\wh{G}_{v_2})) - \dim \centralizer_\lieGroup(\rhoG(\wh{G}_{\eta_1})) \\
    &= \dim \lieGroup + \frac{q-1}{4} - \frac{(q-1)^2}{8} \\
    &<\dim \lieGroup.
   \end{align*}
   Note that the image of $\overline{Y_0}$ is the image of the map $M/H\to \lieGroup$ given by
   $(\tau_{\eta_0},\tau_{\eta_{2}},\tau_{\eta_3})\mapsto\tau_{\eta_0}^{-1}\tau_{\eta_{2}}^{-1}\tau_{\eta_3}^{-1}\rhoG(g_{\eta_3})$.
   Since this map is differentiable we conclude that $\overline{Y}_0$ is not surjective and therefore has degree $0$.
  \end{proof}
 \end{proposition}

  \begin{proposition}\label{degreeY05mod24}
  If $q\equiv 5\pmod {24}$ and $q>5$ then the map $\overline{Y}_0\colon \MMbar_0\to \lieGroup$ has degree $0$.
  \begin{proof}
   By \cref{LowerBoundForDimensionOfIntersectionOfCentralizersCyclicSubgroups} there are matrices
   $A_{\eta_1},A_{\eta_3}\in \lieGroup$ such that $A_{\eta_3}^{-1}A_{\eta_1}$ commutes with $\centralizer_\lieGroup( \rhoG( \wh{G}_{v_2} ))$ and 
   such that the dimension of
   \[K=\left(A_{\eta_3}\centralizer_{\lieGroup}(\rhoG(\wh{G}_{\eta_3}))A_{\eta_3}^{-1}\right)\cap \left(A_{\eta_1}\centralizer_{\lieGroup}(\rhoG(\wh{G}_{\eta_1}))A_{\eta_1}^{-1}\right)\]
   is at least $\frac{m^2}{6\cdot 4}=\frac{(q-1)^2}{96}$.
   Consider the $\HH$-equivariant map $Z\colon \MM_0 \to \lieGroup$ defined by 
   \[\tau\mapsto A_{\eta_3}\tau_{\eta_0}^{-1}\tau_{\eta_{2}}^{-1}\tau_{\eta_{3}}^{-1}A_{\eta_3}^{-1}\cdot A_{\eta_1} \tau_{\eta_1}^{-1}A_{\eta_1}^{-1}\rhoG(g_{\eta_1}).\]
   By \cite[Proposition 5.11]{part1}, the induced maps $\overline{Y}_0,\overline{Z}\colon \MMbar_0 \to \lieGroup$ are homotopic. To conclude, we will prove that $Z$ is not surjective.
   Let
   \begin{align*}
    M &=\left( A_{\eta_3} \centralizer_{\lieGroup}(\rhoG( \wh{G}_{\eta_0})) A_{\eta_3}^{-1}\right) \times \left(A_{\eta_1} \centralizer_{\lieGroup}(\rhoG( \wh{G}_{\eta_1})) A_{\eta_1}^{-1}\right) \\
    &\qquad \times \left(A_{\eta_3} \centralizer_{\lieGroup}(\rhoG( \wh{G}_{\eta_2})) A_{\eta_3}^{-1}\right) \times \left(A_{\eta_3} \centralizer_{\lieGroup}(\rhoG( \wh{G}_{\eta_3})) A_{\eta_3}^{-1}\right) \\
    H &=  \left(A_{\eta_3} \centralizer_{\lieGroup}(\rhoG( \wh{G}_{v_1})) A_{\eta_3}^{-1}\right) \times \left(A_{\eta_3} \centralizer_{\lieGroup}(\rhoG( \wh{G}_{v_3})) A_{\eta_3}^{-1}\right) \times K
   \end{align*}
and consider the free right action $M\rightAction H$ given by
\[(\tau_0,\tau_1,\tau_2,\tau_3)\cdot (\alpha_1,\alpha_3,\alpha_K) =  (\alpha_1^{-1}\tau_0, \tau_1\alpha_{K}, \alpha_3^{-1}\tau_2\alpha_1, \alpha_K^{-1}\tau_3\alpha_3).\]
Finally, note that the image of $Z$ is the image of the $H$-equivariant map $T\colon M\to \lieGroup$
given by $(\tau_0,\tau_1,\tau_2,\tau_3)\mapsto \tau_0^{-1}\tau_2^{-1}\tau_{3}^{-1}\tau_1^{-1}\rho(g_{\eta_1})$,
which cannot be surjective since we have
\begin{align*}
\dim M/H &= \dim M -\dim H \\
&=\dim \MM_0 - \dim \HH + \dim \centralizer_{\lieGroup}( \rhoG( \wh{G}_{v_2} ) ) - \dim K \\
&\leq \dim \lieGroup +\frac{q-1}{4}-\frac{(q-1)^2}{96}\\
\text{(since $q> 25$)}\qquad &< \dim \lieGroup.
\end{align*}
\end{proof}
\end{proposition}

\bibliographystyle{alpha}
\bibliography{references}

\end{document}